\documentclass[letterpaper, 11pt]{amsart}

\usepackage{amsmath,amsthm,amsfonts,amssymb,amscd}
\usepackage{bbm}
\usepackage{bm}
\usepackage{tikz}
\usepackage{tikz-cd}
\usepackage{appendix}
\usepackage{BOONDOX-calo}
\usepackage{standalone}
\usepackage[hidelinks]{hyperref}
\usetikzlibrary{arrows,chains,matrix,positioning,scopes}
\usepackage[letterpaper, left=2.8cm,right=2.8cm, top=2.8cm, bottom=2.8cm]{geometry}
\usepackage{adjustbox}
\usepackage{color}
\usepackage{enumitem}

\newcommand{\CC}{{\mathbb C}}

\newcommand{\ZZ}{{\mathbb Z}}
\newcommand{\QQ}{{\mathbb Q}}
\newcommand{\oQ}{{\overline{\mathbb Q}}}
\newcommand{\K}{{\mathbb K}}
\newcommand{\F}{{\mathbb F}}

\newcommand{\inEnd}{{\underline{End}}}
\newcommand{\fu}{{\mathfrak{u}}}
\newcommand{\fg}{{\mathfrak{g}}}
\newcommand{\inHom}{{\underline{Hom}}}

\newcommand\dual{\raise0.9ex\hbox{$\scriptscriptstyle\vee$}}

\newcommand{\bT}{{\mathbf{T}}}

\theoremstyle{plain}
\newtheorem{thm}{Theorem} 
\newtheorem{prop}[thm]{Proposition}
\newtheorem{lemma}[thm]{Lemma}  
\newtheorem{cor}[thm]{Corollary}
\numberwithin{thm}{section}

\newtheorem{manualtheoreminner}{Theorem}
\newenvironment{thm'}[1]{%
  \renewcommand\themanualtheoreminner{#1}%
  \manualtheoreminner
}{\endmanualtheoreminner}

\theoremstyle{definition}
\newtheorem{defn}[thm]{Definition}

\theoremstyle{remark}
\newtheorem{rem}[thm]{Remark}

\tikzset{>=stealth}

\makeatletter
\def\@seccntformat#1{%
  \protect\textup{\protect\@secnumfont
    \ifnum\pdfstrcmp{subsection}{#1}=0 \bfseries\fi
    \csname the#1\endcsname
    \protect\@secnumpunct
  }%
}  
\makeatother

\makeatletter
\@namedef{subjclassname@2020}{%
  $2020$ Mathematics Subject Classification}
\makeatother

\begin{document}

\title[A remark on Ext groups for motives with maximal unipotent radicals]{A remark on Ext groups for motives with maximal unipotent radicals}
\author{Payman Eskandari}
\address{Department of Mathematics and Statistics, University of Winnipeg, Winnipeg MB, Canada }
\email{p.eskandari@uwinnipeg.ca}
\subjclass[2020]{19E15, 14C15, 18M25, 11G99.}
\begin{abstract}
Let $\bT$ be a neutral tannakian category over a field of characteristic 0. Let $M$ be an object of $\bT$ with a filtration $0=F_0M\subsetneq F_1M\subsetneq \cdots\subsetneq F_kM=M$, such that each successive quotient $F_iM/F_{i-1}M$ is semisimple. Assume that the unipotent radical of the tannakian fundamental group of $M$ is as large as it is permitted under the constraints imposed by the filtration $(F_\bullet M)$. In this note, we first describe the $Ext^1$ groups in the tannakian subcategory of $\bT$ generated by $M$. We then give two applications for motives, one involving 1-motives and another involving mixed Tate motives, leading to some implications of Grothendieck's period conjecture.
\end{abstract}
\maketitle
\vspace{-.2in}
\section{Statement of the result}\label{sec: statement}
Let $\bT$ be a neutral tannakian category (as in \cite{DM82}) over a field $\F$ of characteristic 0. Let $M$ be an object of $\bT$, equipped with a filtration
\begin{equation}\label{eq8}
0=F_0M\subsetneq F_1M\subsetneq \cdots\subsetneq F_kM = M
\end{equation}
such that the associated graded
\[
Gr^F(M)  = \bigoplus\limits_{i=1}^k \, Gr^F_i(M) = \bigoplus\limits_{i=1}^k \, \frac{F_iM}{F_{i-1}M}
\]
is semisimple. The prototype examples to keep in mind are when $\bT$ is the category of graded-polarizable rational mixed Hodge structures or any reasonable tannakian category of mixed motives over a subfield of $\CC$, and $(F_iM)_{0\leq i\leq k}$ is the weight filtration on $M$ with re-labelled indices (with $k$ = the number of weights of $M$). But what we will say below does not require a functorial filtration. Note that even in motivic applications, we do not always have an interesting functorial weight filtration in the categories of interest. For example, no such filtration is known for the de Rham-Betti realization (as in \cite[\S 7.1.6]{An04}).
\medskip\par 
Given a fiber functor $\omega$ for $\bT$ over $\F$, let $\mathcal{U}(M,\omega)$ be the kernel of the natural surjection from the tannakian group of $M$ with respect to $\omega$ (i.e., $\underline{\text{Aut}}^{\otimes}(\omega |_{\langle M\rangle^{\otimes}})$ in the standard notation \cite{DM82}) to the tannakian group of $Gr^F(M)$ with respect to $\omega$. Since $Gr^F(M)$ is semisimple, $\mathcal{U}(M,\omega)$ is the unipotent radical of the tannakian group of $M$ with respect to $\omega$. Tannakian formalism provides us with a canonical subobject 
\[
\fu(M) \subset \inEnd(M) := \inHom(M,M)
\]
(where $\inHom$ means the internal Hom) such that for every choice of a fiber functor $\omega$, the Lie subalgebra
\[
\omega \fu(M) \subset \omega\inEnd(M) = End(\omega M)
\]
is the Lie algebra of $\mathcal{U}(M,\omega)$ (see for instance, \cite[\S 2]{EM2} for a detailed explanation of this). Moreover, since the elements of $\mathcal{U}(M,\omega)$ act trivially on $\omega Gr^F(M)$, we have
\[
\fu(M)\subset F_{-1}\inEnd(M),
\]
where $F_{-1}\inEnd(M)$ is the subobject of $\inEnd(M)$ such that for every fiber functor $\omega$,
\begin{equation}\label{eq9}
\omega F_{-1}\inEnd(M) = \{f\in End_\F(\omega M): f(\omega F_iM)\subset \omega F_{i-1}M \ \ \text{for all $i$}\}.
\end{equation}
(Here and in what follows, for convenience, we set $F_iM= 0$ for $i<0$.)

\begin{defn}\label{def: maximal}
We say $\fu(M)$ or the unipotent radical of the tannakian group of $M$ is {\it maximal} (with respect to the given filtration $(F_\bullet M)$ on $M$) if
\[
\fu(M) = F_{-1}\inEnd(M).
\]
\end{defn}

Some criteria for maximality of $\fu(M)$ in various special cases are given in \cite{BP1}, \cite{Ber01}, \cite{EM2}, \cite{Es1}, and \cite{Es2}. A necessary and sufficient criterion for maximality of $\fu(M)$ in the generality of a filtration $(F_\bullet M)$ in a neutral tannakian category over a field of characteristic zero with no other hypotheses is given in \cite{Es3}, generalizing the criteria given in the previously mentioned references. From arithmetic and algebro-geometric viewpoints, for example, one may be interested in the notion of maximality of $\fu(M)$ in relation to (1) periods and Grothendieck's period conjecture (which predicts the dimension of the motivic Galois group of a motive over $\overline{\QQ}$ should be equal to the transcendence degree of the field generated by the periods of the motive), and (2) the motivic versions of Hodge and Tate conjectures (which predict equality of motivic Galois groups of motives over suitable fields with their Hodge-theoretic and $\ell$-adic analogues, see \cite{An04}). In relation to (1), among all motives with the same associated graded with respect to the weight filtration, the field generated by the periods of a motive with a maximal unipotent radical of motivic Galois group should have the highest transcendence degree. In relation to (2), for motives whose Hodge (resp. $\ell$-adic) realization satisfies this maximality property, the motivic version of the Hodge (resp. Tate) conjecture will hold at the unipotent level for trivial reasons.
\medskip\par 
The purpose of this short note is to make an observation about Yoneda $Ext^1$ groups in tannakian subcategories generated by objects with maximal $\fu$ and draw some consequences of it in relation to Grothendieck's period conjecture. The observation is the following:
\begin{thm}\label{thm1}
Suppose $\bT$ and $M$ are as above (i.e., $\bT$ is a neutral tannakian category over a field $\F$ of characteristic 0 and $M$ is an object of $\bT$ equipped with a finite filtration \eqref{eq8} such that $Gr^FM$ is semisimple). Set
\begin{equation}\label{eq3}
V(M) \ := \ \bigoplus_{i=1}^{k-1}\inHom(Gr^F_{i+1}M, Gr^F_iM).
\end{equation}
If $\fu(M)$ is maximal, then for every semisimple object $X$ in $\langle M\rangle^{\otimes}$ there is a natural (i.e., functorial in $X$) isomorphism
\begin{equation}\label{eq1}
Ext^1_{\langle M\rangle^{\otimes}}(\mathbbm{1}, X) \ \cong \ Hom(V(M), X).
\end{equation}
\end{thm}
Here, $\langle M\rangle^{\otimes}$ means the tannakian subcategory of $\bT$ generated by $M$, and $Ext^1_{\langle M\rangle^{\otimes}}$ means the Yoneda $Ext^1$ group in $\langle M\rangle^{\otimes}$. The notation $Hom$ without any subscript refers to the Hom group in $\bT$ (or $\langle M\rangle^{\otimes}$).

Below we will first consider some special cases of Theorem \ref{thm1}, and then prove the result. Finally, we will give two applications in motivic contexts.

\section{Special cases}
\numberwithin{thm}{subsection}
Unless otherwise indicated, $\bT$ and $M$ are as in Theorem \ref{thm1}. 

\subsection{} Theorem \ref{thm1} has the following implication:
\begin{cor}\label{cor1}
Suppose $\fu(M)$ is maximal. For any $i,j$, the Ext group
\[
Ext^1_{\langle M\rangle^{\otimes}}(Gr^F_j M, Gr^F_i M)
\]
vanishes if and only if $\inHom(Gr^F_jM, Gr^F_iM)$ and $V(M)$ (defined in Theorem \ref{thm1}) have no nonzero isomorphic subobjects.
\end{cor}

This is immediate from Theorem \ref{thm1} in view of semisimplicity of $Gr^F(M)$ and the canonical isomorphism
\[
Ext^1(A, B) \ \cong \  Ext^1(\mathbbm{1},\inHom(A, B))
\]
in a tannakian category.

\subsection{}\label{subsection: filtration is by weights when available} Suppose $\bT$ is filtered, i.e., there exists a filtration $W_\bullet$ (referred to as the weight filtration) on objects of $\bT$ satisfying similar properties to the weight filtration on the category of rational mixed Hodge structures: indexed by $\ZZ$, increasing and finite on every object, functorial, exact, and preserving the tensor structure. By a weight of an object $M$ we mean an integer $p$ such that $W_{p-1}M\neq W_pM$. Let $M$ be an object of $\bT$ with $k$ weights $p_1<\cdots<p_k$ such that its associated graded $Gr^WM$ with respect to the weight filtration is semisimple. In this case, unless otherwise indicated, we will always take the earlier filtration $(F_\bullet M)$ to be the weight filtration with re-labelled indices (i.e., with $F_iM= W_{p_i}M$). Then the object $\fu(M)$ is maximal if and only if $\fu(M)= W_{-1}\inEnd(M)$. The object $V(M)$ defined in Theorem \ref{thm1} is then
\[
\bigoplus\limits_{i=1}^{k-1} \inHom(Gr^W_{p_{i+1}}M, Gr^W_{p_{i}}M).
\]

\subsection{} In the situation of the previous remark, suppose that $M$ has $k=3$ weights $p_1<p_2<p_3$. A comparison of the weights shows that $\inHom(Gr^W_{p_3}M,Gr^W_{p_1}M)$ and 
\[
V(M) = \inHom(Gr^W_{p_2}M,Gr^W_{p_1}M)\oplus \inHom(Gr^W_{p_3}M,Gr^W_{p_2}M)
\]
have no nonzero isomorphic subobjects. We thus obtain:
\begin{cor}\label{cor2}
Suppose $\bT$ is a filtered tannakian category over a field of characteristic 0, $M$ is an object of $\bT$ with 3 weights $p_1<p_2<p_3$ such that $Gr^WM$ is semisimple and $\fu(M)$ is maximal. Then 
\[
Ext^1_{\langle M\rangle^{\otimes}}(Gr^W_{p_3} M, Gr^W_{p_1} M)=0.
\]
\end{cor}

\section{Proof of Theorem \ref{thm1}}
\numberwithin{thm}{section}
Theorem \ref{thm1} follows from the following two lemmas:

\begin{lemma}\label{lem: characterization of Ext groups for algebraic groups}
Let $\mathcal{G}$ be a linear algebraic group over $\F$ (a field of characteristic zero). Denote the category of finite-dimensional representations of $\mathcal{G}$ over $\F$ by $\mathbf{Rep}(\mathcal{G})$. Denote the $Hom$ and Yoneda $Ext^1$ groups for $\mathbf{Rep}(\mathcal{G})$ by $Hom_{\mathcal{G}}$ and $Ext^1_{\mathcal{G}}$. Let $\mathcal{U}$ be the unipotent radical of $\mathcal{G}$. Let $\fu$ be the Lie algebra of $\mathcal{U}$, considered as an object of $\mathbf{Rep}(\mathcal{G})$ through the adjoint action, and let $\fu^{ab}:=\fu/[\fu,\fu]$ be the abelianization of $\fu$, also considered as an object of $\mathbf{Rep}(\mathcal{G})$ via the induced action (note that $\fu^{ab}$ is semisimple, as the action of $\mathcal{G}$ on $\fu^{ab}$ factors through an action of the reductive group $\mathcal{G}/\mathcal{U}$). Then for every semisimple object $X$ of $\mathbf{Rep}(\mathcal{G})$, there is a natural isomorphism
\[
Ext^1_{\mathcal{G}}(\mathbbm{1}, X) \ \cong \ Hom_{\mathcal{G}}(\fu^{ab}, X).
\]
\end{lemma}
\begin{proof}
This result should be well known (see for instance, \cite[\S 16]{Hain23}), although it is not easy to reference to a place in the literature where it is explicitly written in this format. A very explicit proof working with ``bare hands" can be found in \cite{Es2} (see Proposition 3.1.1 therein). For the sake of completeness, here we will include a compact argument using group cohomology for algebraic groups. There are natural isomorphisms
\[
Ext^1_{\mathcal{G}}(\mathbbm{1}, X) \stackrel{(\ast)}{\cong} H^1(\mathcal{G}, X) \stackrel{(\ast\ast)}{\cong} H^1(\mathcal{U}, X)^{\mathcal{G}/\mathcal{U}},
\]
where $(\ast)$ is by the definition of group cohomology and the fact that Yoneda Ext groups agree with those defined using injective resolutions. The isomorphism $(\ast\ast)$ is by the inflation-restriction exact sequence for the $\mathcal{G}$-module $X$ and the subgroup $\mathcal{U}$ of $\mathcal{G}$, on noting that $\mathcal{G}/\mathcal{U}$ is reductive and so $H^i(\mathcal{G}/\mathcal{U}, -)$ vanishes for $i>0$. Since $X$ is a semisimple $\mathcal{G}$-representation, the action of $\mathcal{U}$ on $X$ is trivial, so that
\[H^1(\mathcal{U}, X) \cong Hom_{\mathbf{AlgGr}}(\mathcal{U}, X) \cong Hom_{\mathbf{AlgGr}}(\mathcal{U}^{ab}, X),\]
where $Hom_{\mathbf{AlgGr}}$ is Hom for the category of algebraic groups over $\F$ and $\mathcal{U}^{ab}$ is the abelianization of $\mathcal{U}$. Since $\mathcal{U}^{ab}$ and $X$ are both additive groups, by passing to the Lie algebras we have
\[
Hom_{\mathbf{AlgGr}}(\mathcal{U}^{ab}, X) \cong Hom_\F(\fu^{ab}, X).
\]
The condition of ${\mathcal{G}/\mathcal{U}}$-invariance on $H^1(\mathcal{U}, X)$ translates to ${\mathcal{G}/\mathcal{U}}$-invariance on $Hom_\F(\fu^{ab}, X)$, so that
\[
H^1(\mathcal{U}, X)^{\mathcal{G}/\mathcal{U}} \ \cong \ Hom_{\mathcal{G}/\mathcal{U}}(\fu^{ab}, X) \ = \ Hom_{\mathcal{G}}(\fu^{ab}, X).
\]
\end{proof}

\begin{lemma}\label{lem: abelianization of F_-1End(M)}
Let $\bT$ be a neutral tannakian category over a field $\F$ of characteristic 0 and $M$ an object of $\bT$ equipped with a finite filtration \eqref{eq8} such that $Gr^FM$ is semisimple. Then
\[
\bigm(F_{-1}\inEnd(M)\bigm)^{ab} \ \cong \ V(M).
\]
(See equations \eqref{eq9} and \eqref{eq3} for the definitions of $F_{-1}\inEnd(M)$ and $V(M)$. Here and in what follows, for a Lie algebra object $\fg$ of $\bT$, the abelianization $\fg/[\fg,\fg]$ (also an object of $\bT$) is denoted by $\fg^{ab}$.)
\end{lemma}
\begin{proof}
Let $F_{-2}\inEnd(M)$ be the subobject of $\inEnd(M)$ whose image under any fiber functor $\omega$ consists of the linear endomorphisms $f$ of $\omega M$ such that $f(\omega F_iM)\subset \omega F_{i-2}M$ for every $i$ (where as before, we set $F_i M = 0$ for $i<0$). It is clear from definitions that $[F_{-1}\inEnd(M),F_{-1}\inEnd(M)]\subset F_{-2}\inEnd(M)$. We shall show that $F_{-2}\inEnd(M)=[F_{-1}\inEnd(M),F_{-1}\inEnd(M)]$ and construct a surjective morphism
\begin{equation}\label{eq10}
F_{-1}\inEnd(M) \ \longrightarrow \ \bigoplus_{i=1}^{k-1} \, \inHom(Gr^F_{i+1}M, Gr^F_iM)
\end{equation}
whose kernel is $F_{-2}\inEnd(M)$.

Choose a fiber functor $\omega$ over $\F$. Set $F_i\omega M:=\omega F_iM$ for each $i$. Let $\mathcal{G}(M,\omega)$ be the tannakian group of $M$ with respect to $\omega$. For each $2\leq i\leq k$ choose a linear section $s_i$ of the map $F_i \omega M\twoheadrightarrow Gr^F_i\omega M$ to identify $F_i \omega M\cong F_{i-1} \omega M \oplus Gr^F_i\omega M$ as vector spaces. By successive applications of these we have identifications
\[
\omega M \ = \ F_k \omega M \ \cong \ F_{k-1}\omega M\oplus Gr^F_k\omega M \ \cong \ F_{k-2}\omega M\oplus Gr^F_{k-1}\omega M \oplus Gr^F_k\omega M \ \cong \ \cdots 
\]
and finally an identification
\begin{equation}\label{eq4}
\omega M  \ \cong \ \bigoplus\limits_{i=1}^k Gr^F_i\omega M
\end{equation}
and subsequently an identification
\begin{equation}\label{eq5}
End_\F(\omega M) \ \stackrel{\eqref{eq4}}{\cong} \ End_\F(\bigoplus\limits_{i=1}^k Gr^F_i\omega M) \ \cong \ \bigoplus\limits_{1\leq i,j\leq k} Hom_\F(Gr^F_j\omega M, Gr^F_i\omega M).
\end{equation}
Via these identifications, we have $F_j \omega M\cong \bigoplus\limits_{i=1}^j Gr^F_i\omega M$ and
\begin{equation}\label{eq6}
\omega F_{-\ell}\inEnd(M) \ \cong \ \bigoplus\limits_{j-i\geq \ell} Hom_\F(Gr^F_j\omega M, Gr^F_i\omega M) \hspace{.5in} (\ell=1,2).
\end{equation}
An easy exercise in linear algebra now verifies that
\[
\omega F_{-2}\inEnd(M) \subset [\omega F_{-1}\inEnd(M), \omega F_{-1}\inEnd(M)],
\]
hence giving the equality of $F_{-2}\inEnd(M)$ and the derived algebra of $F_{-1}\inEnd(M)$. 

Let
\[
\pi: \omega F_{-1}\inEnd(M) \ \twoheadrightarrow \ \bigoplus\limits_{i = 1}^{k-1} Hom_\F(Gr^F_{i+1}\omega M, Gr^F_i\omega M)
\]
be the composition of the identification \eqref{eq6} for $\ell=1$ with the projection to the components with $j=i+1$. From \eqref{eq6} the kernel of $\pi$ is $\omega F_{-2}\inEnd(M)$. We claim that (i) the map $\pi$ is independent of the choices of the sections $(s_i)$ used to obtain identifications \eqref{eq4} and \eqref{eq5}, and (ii) $\pi$ is equivariant under the actions of $\mathcal{G}(M,\omega)$ on $\omega F_{-1}\inEnd(M)$ and $\bigoplus_{i = 1}^{k-1} Hom_\F(Gr^F_{i+1}\omega M, Gr^F_i\omega M)$ corresponding to the objects $F_{-1}\inEnd(M)$ and $\bigoplus_{i = 1}^{k-1} \inHom(Gr^F_{i+1}M, Gr^F_iM)$ under the equivalence of categories 
\begin{equation}\label{eq7}
\langle M\rangle^{\otimes} \ \xrightarrow{ \ \omega \ } \ \mathbf{Rep}(\mathcal{G}(M,\omega)).
\end{equation}
This will complete the proof of the lemma. 

Write elements of $End_\F(Gr^F \omega M)$ as $k\times k$ matrices via the canonical decomposition  
\[
End_\F(Gr^F \omega M) \ \cong \ \bigoplus\limits_{1\leq i,j\leq k} Hom_\F(Gr^F_j\omega M, Gr^F_i\omega M).
\]
Under the isomorphism \eqref{eq5}, $\omega F_{-1}\inEnd(M)$ consists of all the strictly upper-triangular elements of $End_\F(Gr^F \omega M)$. The isomorphisms \eqref{eq5} for two choices of sections $(s_i)$ and $(s'_i)$ differ by conjugation by an upper-triangular element of $End_\F(Gr^F \omega M)$ whose diagonal entries are all identity maps. Such a conjugation leaves the super-diagonal entries $f_{i,i+1}$ of a strictly upper-triangular matrix $(f_{i,j})$ unchanged. This shows the independence of $\pi$ from our choice of sections. 

As for $\mathcal{G}(M,\omega)$-equivariance, the action of $\mathcal{G}(M,\omega)\subset GL(\omega M)$ on $\omega F_{-1}\inEnd(M)$ (resp. $Hom_\F(Gr^F_j\omega M, Gr^F_i\omega M)$) corresponding to the object $F_{-1}\inEnd(M)$ (resp. $\inHom(Gr_j M, Gr_i M)$) under \eqref{eq7} is given by conjugation (resp. by $\sigma\cdot f_{i,j} = \sigma_{Gr^F_i M}\circ f_{i,j} \circ {\sigma_{Gr^F_j M}}^{-1}$, where $\sigma_X$ denotes the action of $\sigma$ on $\omega X$ for any object $X$ of $\langle M\rangle^{\otimes}$). The subgroup $\mathcal{G}(M,\omega)$ of $GL(\omega M)$ respects the filtration $F_\bullet$ on $\omega M$. Identifying $GL(\omega M)$ with $GL(Gr^F \omega M)$ via \eqref{eq4} (for a choice of sections $(s_i)$) and writing its elements as $k\times k$ matrices, $\mathcal{G}(M,\omega)$ is thus contained in the upper-triangular subgroup of $GL(Gr^F \omega M)$. Given $\sigma=(\sigma_{i,j}) \in \mathcal{G}(M,\omega)$ and $f=(f_{i,j})\in \omega F_{-1}\inEnd(M)$ (with $\sigma$ upper-triangular and $f$ strictly so), a direct computation shows that the super-diagonal entry $(i,i+1)$ of $\sigma f \sigma^{-1}$ is $\sigma_{i,i}\circ f_{i,i+1}\circ \sigma_{i+1,i+1}^{-1}$. On noting that $\sigma_{Gr^F_j M}=\sigma_{j,j}$ for all $j$, we thus have $\pi(\sigma\cdot f) = \sigma\cdot(\pi(f))$, as desired. 
\end{proof}

We are ready to deduce Theorem \ref{thm1}. Let $\bT$ and $M$ be as in the statement of the theorem (with $\fu(M)$ not yet assumed to be maximal). Let $\omega$ a fiber functor for $\bT$ over $\F$. Apply Lemma \ref{lem: characterization of Ext groups for algebraic groups} taking $\mathcal{G}$ to be the tannakian group $\mathcal{G}(M,\omega)$ of $M$ with respect to $\omega$. The unipotent radical of $\mathcal{G}$ is the group $\mathcal{U}(M,\omega)$ with Lie algebra $\omega\fu(M)$. The representations of $\mathcal{G}$ on $\omega\fu(M)$ and $\omega\fu(M)^{ab}$ under the equivalence of categories \eqref{eq7} are induced by the adjoint action of $\mathcal{G}$. In view of the same equivalence of categories and Lemma \ref{lem: characterization of Ext groups for algebraic groups}, for every semisimple object $X$ of $\langle M\rangle^{\otimes}$  we thus have a natural isomorphism
\begin{equation}\label{eq2}
Ext^1_{\langle M\rangle^{\otimes}}(\mathbbm{1}, X) \ \cong \ Hom(\fu(M)^{ab}, X).
\end{equation}
When $\fu(M)$ is maximal, combining with Lemma \ref{lem: abelianization of F_-1End(M)} we obtain Theorem \ref{thm1}.

\begin{rem}
The construction of the isomorphism \eqref{eq1} of Theorem \ref{thm1} given above depends on the choice of a fiber functor $\omega$ for $\bT$. The isomorphism itself however is independent of this choice. Indeed, a detailed proof of the fact that the isomorphism \eqref{eq2} is independent from the choice of $\omega$ can be found in \cite{Es2} (see Lemma 3.2.1 therein). That the isomorphism of Lemma \ref{lem: abelianization of F_-1End(M)} is independent of the choice of $\omega$ can be seen from a more natural description of the map \eqref{eq10} along the lines of \cite[\S 2.6]{Es2} (where this more natural description is given when $k=3$). 
\end{rem}

\section{Applications}
In what follows, we will give two examples related to motives. Given a tannakian category $\mathbf{MM}(\K)$ of mixed motives over a field $\K\subset\CC$, by the motivic Galois group of an object $M$ of $\mathbf{MM}(\K)$ we mean the tannakian group of $M$ with respect to the fiber functor given by the Betti realization (with coefficients in $\QQ$). If $\K$ is algebraic over $\QQ$, by Grothendieck's period conjecture (in the context of the category $\mathbf{MM}(\K)$) we mean the following statement: for every object $M$ of $\mathbf{MM}(\K)$, the dimension of the motivic Galois group of $M$ is equal to the transcendence degree of the field generated over $\K$ (or $\QQ$) by the periods of $M$. We refer to this statement for a fixed object $M$ of $\mathbf{MM}(\K)$ as Grothendieck's period conjecture for $M$. The reader can consult an appendix of \cite{Be20} due to Andr\'{e} for an excellent exposition of Grothendieck's conjecture. Note that the statement of the conjecture involves the choice of a tannakian category of motives.

Denote the category of rational mixed Hodge structures by $\mathbf{MHS}$. The Mumford-Tate group of a rational mixed Hodge structure $H$ is the tannakian group of $H$ in $\mathbf{MHS}$ with respect to the forgetful fiber functor. Any reasonable candidate for $\mathbf{MM}(\K)$ for any subfield $\K\subset\CC$ is equipped with a Hodge realization functor $\mathfrak{R}^H:\mathbf{MM}(\K)\rightarrow \mathbf{MHS}$, which is exact faithful linear and respects tensor products. For every motive $M$, the functor $\mathfrak{R}^H$ identifies the Mumford-Tate group of $\mathfrak{R}^HM$ (or for simplicity, of $M$) as a subgroup of the motivic Galois group of $M$.

\subsection{Application I: 1-motives}
Fix an embedding $\oQ\subset\CC$. We will take the background on the theory of Deligne 1-motives \cite{De74} for granted. Let $M=[\ZZ^n\xrightarrow{\mu} G]$ be a 1-motive over $\oQ$, where $G$ is an extension of a nonzero abelian variety $A$ by a torus $\mathbb{G}_m^s$. Denote the dual abelian variety to $A$ by $A^\ast$. Let 
\[
\underline{Q} = (Q_1,\ldots, Q_s)  \in A^\ast(\oQ)^s \hspace{.5in}(\text{resp.} \ \ \underline{P} = (P_1,\ldots, P_n)  \in A(\oQ)^n)
\]
correspond to $G$ (resp. $M/\mathbb{G}_m^s = [\ZZ^n\xrightarrow{\mu} G \rightarrow A]$) under the canonical isomorphism
\[
Ext^1(A, \mathbb{G}_m) \cong A^\ast(\overline{\QQ}) \hspace{.5in}(\text{resp.} \ \ Ext^1(\ZZ, A) \cong A(\overline{\QQ})),
\]
where the Ext groups are in the category of 1-motives over $\oQ$.

Let $\mathbf{NM}(\oQ)$ be Nori's tannakian category of (non-effective) motives over $\overline{\QQ}$ \cite{HM17}. Let $\mathfrak{R}^H: \mathbf{NM}(\oQ)\rightarrow \mathbf{MHS}$ be the Hodge realization functor. We denote the Nori motive associated to our 1-motive $M$ also by the same symbol. By \cite[Corollary 6.7]{ABV15}, $\mathfrak{R}^H M$ (resp. de Rham and $\ell$-adic realizations of $M$ as a Nori motive) is the usual Hodge (resp. de Rham and $\ell$-adic) realization of the 1-motive $M$ as constructed by Deligne in \cite{De74}.

Andr\'{e} \cite{An19} has proved that the motivic Galois group of $M$ (in the context of $\mathbf{NM}(\oQ)$) coincides with its Mumford-Tate group. Moreover, by \cite[Theorem 6.2]{Jo14} (or \cite[Corollary 3.8.1]{Es2}, in view of the fullness of the Hodge realization functor on the category of 1-motives over $\oQ$ up to isogeny \cite[Proposition 2.1]{An19}) the following two statements are equivalent:
\begin{itemize}
\item[(i)] The unipotent radical of the Mumford-Tate group of $M$ is maximal. That is, in the notation of \S \ref{sec: statement}, $\fu(\mathfrak{R}^H M)$ is maximal, where $\bT$ is taken to be $\mathbf{MHS}$ (and the filtration on $\mathfrak{R}^H M$ is the weight filtration with relabelled indices, see \S \ref{subsection: filtration is by weights when available}).
\item[(ii)] The only complex abelian subvariety of $A^n\times {A^\ast}^s$ that contains the point
\[k\cdot (\underline{P},\underline{Q}) \ \in \ A^n\times {A^\ast}^s\]
for some nonzero integer $k$ is $A^n\times {A^\ast}^s$.
\end{itemize}

As an application of Theorem \ref{thm1} one obtains the following:
\begin{prop}\label{prop: algebraic indep of log values from periods of a 1-motive with maximal u}
Suppose that the 1-motive $M$ satisfies the equivalent conditions (i) and (ii) above. Suppose $r\in \oQ^\times$ is not a root of unity. Assuming Grothendieck's period conjecture, (any value of) $\log(r)$ is algebraically independent from the periods of $M$. 
\end{prop}
\begin{proof}
Statement (i) implies that $\fu(M)$ (for the motivic Galois group) is maximal. By Theorem \ref{thm1}, the Ext group $Ext^1(\mathbbm{1},\QQ(1))$ vanishes for the subcategory $\langle M\rangle^{\otimes}$ of $\mathbf{NM}(\oQ)$. Now consider the Kummer motive $L_r:=H_1(\mathbb{G}_m,\{1,r\})$ in $\mathbf{NM}(\oQ)$. It is a nontrivial extension of $\mathbbm{1}$ by $\QQ(1)$ and hence does not belong to $\langle M\rangle^{\otimes}$. Thus the natural surjection from the motivic Galois group of $M' := M \oplus L_r$ to the motivic Galois group of $M$ is not an isomorphism. We have a commutative diagram
\[
\begin{tikzcd}
1 \arrow[r] & \mathcal{U}(M', \omega_B) \arrow[d, twoheadrightarrow] \arrow[r] & \mathcal{G}(M', \omega_B) \arrow[d, twoheadrightarrow] \arrow[r] & \mathcal{G}(Gr^WM', \omega_B) \ar[equal,d] \arrow[r] & 1\\
1 \arrow[r] & \mathcal{U}(M, \omega_B) \arrow[r] & \mathcal{G}(M, \omega_B) \arrow[r] & \mathcal{G}(Gr^WM, \omega_B) \arrow[r] & 1
\end{tikzcd}
\]
where $\omega_B$ is the Betti functor on $\mathbf{NM}(\oQ)$ and $\mathcal{G}(-,\omega_B)$ denotes the motivic Galois group, and the maps are either inclusion or by restriction. It follows that the surjection $\mathcal{U}(M', \omega_B)\rightarrow \mathcal{U}(M, \omega_B)$ is not an isomorphism. Since the groups $\mathcal{U}(-, \omega_B)$ are unipotent, we obtain that $\mathcal{U}(M', \omega_B)$ has a larger dimension than $\mathcal{U}(M, \omega_B)$. Assuming Grothendieck's period conjecture for $M'$, it follows that the transcendence degree of the field generated by the periods of $M'$ must be larger than the transcendence degree of the field generated by the periods of $M$. This gives the claim since $L_r$ has a period matrix with entries $1/(2\pi i)$, $\log(r)/(2\pi i)$, $0$ and $1$. (Note that $1/(2\pi i)$ is already in the algebra of periods of $M$.)
\end{proof}

\begin{rem}
Here we worked with Nori's category to make the discussion concrete. Proposition \ref{prop: algebraic indep of log values from periods of a 1-motive with maximal u} is valid in the setting of any candidate for the tannakian category $\mathbf{MM}(\oQ)$ of mixed motives over $\oQ$, as long as pure objects of $\mathbf{MM}(\oQ)$ are semisimple and the Hodge and de Rham realization functors on 1-motives, considered as objects of $\mathbf{MM}(\oQ)$, agree with those constructed in \cite{De74}. The equality of the motivic Galois and Mumford-Tate groups for arbitrary 1-motives over $\oQ$ (which is specifically known in the context of $\mathbf{NM}(\oQ)$) is not needed by the argument for Proposition \ref{prop: algebraic indep of log values from periods of a 1-motive with maximal u}: the assumption of maximality of $\fu(\mathfrak{R}^HM)$ forces the maximality of $\fu(M)$.
\end{rem}

\begin{rem}
We refer the reader to \cite[\S 5.2]{BP1} and \cite[\S 2]{Be20} for computations of period matrices of 1-motives, particularly in the case when $A$ is an elliptic curve and $n=s=1$. 
\end{rem}

\subsection{Application II: Mixed Tate motives}
What follows here supplements the examples in \cite[\S 6.9]{EM2} and \cite[\S 5.5]{Es1}. In this second application we will work in the category $\mathbf{MTM}(\QQ)$ of mixed Tate motives over $\QQ$, constructed by Levine \cite{Le93} or Voevodsky \cite{Vo00}. As it is well known, thanks to Borel \cite{Bo74}, Beilinson \cite{Bei84} and the relationship between the Ext groups in $\mathbf{MTM}(\QQ)$ and the K-theory of $\QQ$, for every odd integer $n\geq 3$ the vector space $Ext^1(\mathbbm{1},\QQ(n))$ (where $Ext^1$ here is for $\mathbf{MTM}(\QQ)$) is 1-dimensional, generated by an object with a period matrix with entries $1/(2\pi i)^n$, $\zeta(n)/(2\pi i)^n$, $0$, $1$, where $\zeta$ is Riemann's zeta function; we will denote this object by $Z_n$. The vector space $Ext^1(\mathbbm{1},\QQ(1))$ is isomorphic to $\QQ^\times\otimes\QQ$, with the extension corresponding to $r\otimes 1\in \QQ^\times\otimes\QQ$ given by the weight filtration on the Kummer motive $L_r=H_1(\mathbb{G}_m,\{1,r\})$. The groups $Ext^1(\mathbbm{1},\QQ(n))$ for all other values of $n$, as well as all the $Ext^2$ groups vanish.

Fix two distinct odd integers $a,b>1$ and $r\in\QQ^\times\setminus\{1,-1\}$. Since the $Ext^2$ groups in $\mathbf{MTM}(\QQ)$ all vanish, using Grothendieck's formalism of blended extensions ({\it extensions panach\'{e}es} \cite{Gr68}) one can see that there exists an object $M$ of $\mathbf{MTM}(\QQ)$ such that
\[
Gr^WM \ \simeq \ \QQ(1+a+b)\oplus \QQ(1+a)\oplus \QQ(1)\oplus \mathbbm{1} 
\]
and
\begin{equation}\label{eq11}
W_{-2-2a}M \simeq Z_{b}(1+a) , \ \ W_{-2}M/W_{-2-2a-2b}M \simeq Z_a(1), \ \ \text{and} \ \ M/W_{-2-2a}M \simeq L_r
\end{equation}
(see \cite[Theorem 3.3.1(c)]{Es1}). Then $\fu(M)$ is maximal. Indeed, this follows from \cite[Theorem 5.3.2]{Es1} (note that the ``graded-independence" hypothesis therein holds for our $M$ thanks to its weights). Alternatively, the maximality of $\fu(M)$ can also be seen using the isomorphism \eqref{eq2}: The Ext groups $Ext^1_{\langle M\rangle^{\otimes}}(\mathbbm{1},\QQ(n))$ for $n=1,a,b$ do not vanish because of the nonsplit subquotients of $M$ given in \eqref{eq11}, so that $\fu(M)^{ab}$ must have $\QQ(1)$, $\QQ(a)$ and $\QQ(b)$ as subquotients. From this one can deduce that $\fu(M)$ is maximal.
\medskip\par 
It now follows from Corollary \ref{cor1} that 
\[
Ext^1_{\langle M\rangle^{\otimes}}(\mathbbm{1},\QQ(1+a+b)) = 0.
\]
Thus $Z_{1+a+b}$ does not live in the category $\langle M\rangle^{\otimes}$. A similar argument to the one for Proposition \ref{prop: algebraic indep of log values from periods of a 1-motive with maximal u} now shows that Grothendieck's period conjecture for $M\oplus Z_{1+a+b}$ implies that $\zeta(1+a+b)$ is algebraically independent from the periods of $M$.

\begin{rem}
Because of the subquotient $L_r$, the mixed Tate motive $M$ is ramified at $r$ and so the nature of some of its periods is not known, although one has Goncharov's conjectural description \cite[Conjecture 0.17(a)]{Go95}. See \cite[\S 1 \& \S 6.8]{EM2} (also \cite[Remark 5.5.1]{Es1}) and \cite[\S 1.9 \& \S 3]{DW16} for further discussions related to this.
\end{rem}

\begin{rem}
It is worth mentioning that the object $M$ of $\mathbf{MTM}(\QQ)$ satisfying \eqref{eq11} is not unique. In fact, the set of isomorphism classes of such objects of $\mathbf{MTM}(\QQ)$ is a torsor over $Ext^1(\mathbbm{1},\QQ(1+a+b))\simeq \QQ$. This follows from Theorem 3.3.1(b) and Proposition 4.2.3 of \cite{Es1}, in view of the fact that $Ext^1(\mathbbm{1},\QQ(2k))$ vanishes for all $k$ - see \S 5.5. of the same reference for more details.
\end{rem}

\end{document}